\documentclass[10pt]{article}

\usepackage{amsmath}
\usepackage{amsfonts}
\usepackage{amssymb}
\usepackage[english]{babel}
\usepackage{amsthm}

\newcommand{\mR}{\mathbb{R}}
\newcommand{\mC}{\mathbb{C}}
\newcommand{\mN}{\mathbb{N}}

\newcommand{\mZ}{\mathbb{Z}}
\newcommand{\p}{\partial}
\newcommand{\ux}{\underline{x}}

\newcommand{\uom}{\underline{\omega}}
\newcommand{\pux}{\underline{\partial}}

\newcommand{\mcH}{\mathcal{H}}

\newcommand{\enb}{\overline{e_0}}
\newcommand{\Dbar}{\overline{D}}
\newcommand{\pxn}{\partial_{x_0}}
\newcommand{\onehalf}{\frac{1}{2}}
\newcommand{\la}{\langle}
\newcommand{\ra}{\rangle}

\newcommand{\dps}{\displaystyle}
\newcommand{\dpsf}{\displaystyle\frac}

\newtheorem{lemma}{Lemma}
\newtheorem{proposition}{Proposition}
\newtheorem{remark}{Remark}

\numberwithin{equation}{section}
\numberwithin{theorem}{section}
\numberwithin{proposition}{section}
\numberwithin{lemma}{section}
\numberwithin{definition}{section}
\numberwithin{remark}{section}
\numberwithin{corollary}{section}

\begin{document}

\title{Representation of Distributions by Harmonic and Monogenic Potentials in Euclidean Space}

\author{F.\ Brackx, H.\ De Bie, H.\ De Schepper}

\date{\small{Clifford Research Group, Department of Mathematical Analysis,\\ Faculty of Engineering and Architecture, Ghent University\\
Building S22, Galglaan 2, B-9000 Gent, Belgium\\}}

\maketitle

\begin{abstract}
\noindent In the framework of Clifford analysis, a chain of harmonic and monogenic potentials in the upper half of Euclidean space $\mR^{m+1}_+$ was recently constructed, including a higher dimensional analogue of the logarithmic function in the complex plane, and their distributional boundary values were computed. In this paper we determine those potentials in lower half--space $\mR^{m+1}_-$ and investigate whether they can be extended through the boundary $\mR^m$. This is a stepstone to the representation of a doubly infinite sequence of distributions in $\mR^m$, consisting of positive and negative integer powers of the Dirac and the Hilbert--Dirac operator, as the {\em jump} across $\mR^m$ of monogenic functions in the upper and lower half--spaces, in this way providing a sequence of interesting examples of Clifford hyperfunctions.
\end{abstract}


\section{Introduction}
\label{intro}


Hyperfunctions are localizable generalized functions; they form a generalization of the notion of distribution. Their history goes back to the works of G.\ K\"othe (\cite{koethe}), H.G.\ Tillman (\cite{tillman}), et al. and culminated from the 1960's on in the works of the Japanese school including M.\  Sato (\cite{sato}), H.\ Komatsu (\cite{komatsu}), M.\  Morimoto (\cite{morimoto}), et al. One of the construction methods for  a hyperfunction on the real line is to consider the boundary values of a holomorphic function in both the upper and lower complex half--planes, the hyperfunction itself then being the equivalence class of the difference of this holomorphic function across the real axis. Typical examples of one--dimensional hyperfunctions are the Heaviside function $Y(x)$, the delta or Dirac distribution $\delta(x)$,  and the Hilbert kernel or Cauchy principal value distribution $H(x) = - \frac{1}{\pi}{\rm Pv} \frac{1}{x}$, showing the following hyperfunction representations (the branching line for the logarithmic function being taken on the negative real axis):
\begin{equation}
\label{hyperheaviside}
Y(x) \longleftrightarrow (- \frac{1}{2\pi i} \ln{(-z)}, - \frac{1}{2\pi i} \ln{(-z)})
\end{equation}
\vspace*{1mm}
\begin{equation}
\label{hyperdelta}
\delta(x) \longleftrightarrow (- \frac{1}{2\pi i} \frac{1}{z}, - \frac{1}{2\pi i} \frac{1}{z})
\end{equation}
\vspace*{1mm}
\begin{equation}
\label{hyperhilbert}
- \frac{1}{\pi}{\rm Pv} \frac{1}{x}  \longleftrightarrow (- \frac{1}{2\pi} \frac{1}{z},  \frac{1}{2\pi} \frac{1}{z})
\end{equation}
\vspace*{1mm}

\noindent
The case of several variables was developed by M.\ Sato \cite{sato} using cohomology theory. In \cite{sommen1, sommen2} F.\ Sommen established a valuable and elegant alternative theory of multidimensional hyperfunctions within the context of Clifford analysis. Clifford analysis has become an independent discipline of classical analysis; roughly speaking it is a function theory for functions defined in Euclidean space $\mR^{m+1}$ and taking their values in (subspaces of) the universal Clifford algebra $\mR_{0,m+1}$ constructed over $\mR^{m+1}$, equipped with a quadratic form of signature $(0,m+1)$. The concept of a higher dimensional holomorphic function, mostly called {\em monogenic} function, is expressed by means of a generalized Cauchy--Riemann operator, which is a combination of the derivative with respect to one of the real variables, say $x_0$, and the so--called Dirac operator $\pux$ in the remaining real variables $(x_1, x_2, \ldots, x_m)$. The generalized Cauchy--Riemann operator $D$ and its Clifford algebra conjugate $\overline{D}$ factorize the Laplace operator, whence Clifford analysis may be seen as a refinement of harmonic analysis.\\

In a recent paper \cite{bdbds1} a generalization was constructed of the  logarithmic function $\ln{z}$ to Euclidean upper half--space $\mR^{m+1}_+$. The construction of this higher dimensional monogenic logarithmic function was carried out in the framework of Clifford analysis, its starting point being the fundamental solution of the aforementioned generalized Cauchy--Riemann operator $D$, also called Cauchy kernel, and its relation to the Poisson kernel and its harmonic conjugate in $\mR^{m+1}_+$. We then proceeded by induction in two directions, {\em downstream} by differentiation and {\em upstream} by primitivation, yielding a doubly infinite chain of monogenic, and thus harmonic, potentials. This chain mimics the well--known sequence of holomorphic potentials in $\mC_+$ (see e.g. \cite{slang}):
$$
\frac{1}{k!} z^k \left[ \ln z - ( 1 + \frac{1}{2} + \ldots + \frac{1}{k}) \right] \rightarrow \ldots \rightarrow z ( \ln z - 1) \rightarrow \ln z 
\stackrel{\frac{d}{dz}}{\longrightarrow} \frac{1}{z} \rightarrow - \frac{1}{z^2} \rightarrow \ldots \rightarrow (-1)^{k-1} \frac{(k-1)!}{z^k}
$$
Identifying the boundary of upper half--space with $\mR^m \cong \{(x_0,\ux) \in \mR^{m+1} : x_0 = 0\}$, the distributional limits for $x_0 \rightarrow 0+$ of those potentials were computed. They split up into two classes of distributions, which are linked by the Hilbert transform, one scalar--valued, the second one Clifford vector--valued. They belong to two of the four families of Clifford distributions which were thoroughly studied in a series of papers, see \cite{fb1, fb2, distrib} and the references therein. More particularly half of them may be recovered as fundamental solutions of specific powers of the Dirac operator, and also half of them, however not the complementary ones, as fundamental solutions of specific powers of the Laplace operator. By introducing two new pseudodifferential operators, next to and related to the complex powers of the Dirac and Laplace operators, the whole doubly infinite set of distributional boundary values may now be identified as fundamental solutions of the four operators. As a remarkable demonstration of symmetry, the distributional boundary values also can serve as convolution kernels for the corresponding pseudodifferential operators of the same kind but with opposite exponent.  This boundary behavior of the harmonic and monogenic potentials was studied, however restricted to approaching $\mR^m$ from upper half--space, in \cite{bdbds2}.\\

In this paper we complete the study of the boundary behaviour of the potentials by considering distributional limits approaching the boundary $\mR^m$ from lower half--space $\mR^{m+1}_-$. This enables us to express the doubly infinite sequence of distributional boundary values in $\mR^m$ as hyperfunctions involving the aforementioned monogenic potentials. In particular we obtain the multidimensional analogues of the one--dimensional hyperfunctions (\ref{hyperheaviside}), (\ref{hyperdelta}) and (\ref{hyperhilbert}). We say that our hyperfunction representations are direct as each of the considered distributions is linked to one specific monogenic potential and we have not to recur to the standard Cauchy transform by means of which distributions can be represented as hyperfunctions. Remarkably, the parity of the dimension $m$ plays a crucial role in these direct representations. If $m$ is even, the direct hyperfunction representation involving the upstream potentials is lost, leaving the Cauchy representation as the only alternative in this case.\\

The organization of the paper is as follows. To make the paper self--contained we recall in Section 2 the basics of Clifford algebra and Clifford analysis. In Section 3 we study the boundary behaviour of the harmonic and monogenic potentials when approaching the boundary $\mR^m$ from the lower half--space. Then we express each of the obtained boundary distributions as a hyperfunction involving downstream monogenic potentials (Section 4) and upstream monogenic potentials (Section 5).


\section{Prerequisites of Clifford analysis}
\label{basics}


Clifford analysis (see e.g. \cite{red,dss,gilmur,guerli}) is a function theory which offers a natural and elegant generalization to higher dimension of holomorphic functions in the complex plane and refines harmonic analysis. Let  $(e_0, e_1,\ldots,e_m)$ be the canonical orthonormal basis of Euclidean space $\mR^{m+1}$ equipped with a quadratic form of signature $(0,m+1)$. Then the non--commutative multiplication in the universal real Clifford algebra $\mR_{0,m+1}$ is governed by the rule 
$$
e_{\alpha} e_{\beta} + e_{\beta} e_{\alpha} = -2 \delta_{\alpha \beta}, \qquad \alpha,\beta = 0, 1,\ldots,m
$$
whence $\mR_{0,m+1}$ is generated additively by the elements $e_A = e_{j_1} \ldots e_{j_h}$, where $A=\lbrace j_1,\ldots,j_h \rbrace \subset \lbrace 0,\ldots,m \rbrace$, with $0\leq j_1<j_2<\cdots < j_h \leq m$, and $e_{\emptyset}=1$. 
For an account on Clifford algebra we refer to e.g. \cite{porteous}.\\[-2mm]

We identify the point $(x_0, x_1, \ldots, x_m) \in \mR^{m+1}$ with the Clifford--vector variable 
$$
x = x_0 e_0  + x_1 e_1  + \cdots x_m e_m = x_0 e_0  + \ux
$$ 
and the point $(x_1, \ldots, x_m) \in \mR^{m}$ with the Clifford--vector variable $\ux$. 
The introduction of spherical co--ordinates $\ux = r \uom$, $r = |\ux|$, $\uom \in S^{m-1}$, gives rise to the Clifford--vector valued locally integrable function $\uom$, which is to be seen as the higher dimensional analogue of the {\em signum}--distribution on the real line; we will encounter $\uom$ as one of the distributions discussed below.\\[-2mm]

At the heart of Clifford analysis lies the so--called Dirac operator 
$$
\p = \pxn e_0 + \p_{x_1} e_1 + \cdots \p_{x_m} e_m =  \pxn e_0 + \pux
$$
which squares to the negative Laplace operator: $\p^2 = - \Delta_{m+1}$, while also $\pux^2 = - \Delta_{m}$. The fundamental solution of the Dirac operator $\p$ is given by
$$
E_{m+1} (x) = - \frac{1}{\sigma_{m+1}} \ \frac{x}{|x|^{m+1}}
$$
where $\sigma_{m+1} = \frac{2\pi^{\frac{m+1}{2}}}{\Gamma(\frac{m+1}{2})}$ stands for the area of the unit sphere $S^{m}$ in $\mR^{m+1}$.
We also introduce the generalized Cauchy--Riemann operator 
$$
D = \onehalf \enb \p = \onehalf (\pxn + \enb \pux)
$$ 
which, together with its Clifford algebra conjugate $\Dbar = \onehalf(\pxn - \enb \pux)$, also decomposes the Laplace operator: $D \Dbar = \Dbar D = \frac{1}{4} \Delta_{m+1}$. \\[-2mm]

A continuously differentiable function $F(x)$, defined in an open region $\Omega \subset \mR^{m+1}$ and taking values in the Clifford algebra  $\mR_{0,m+1}$, is called (left--)monogenic if it satisfies  in $\Omega$ the equation $D F = 0$, which is equivalent with $\p F = 0$.\\

We will extensively use two families of distributions in $\mR^m$, which have been thoroughly studied in \cite{fb1,fb2,distrib}. The first family $\mathcal{T} = \{ T_\lambda : \lambda \in \mC\}$ is very classical (see e.g. \cite{ries, helgason}). It consists of the radial distributions 
$$
T_\lambda = {\rm Fp} \ r^{\lambda} = {\rm Fp} \ (x_1^2 + \ldots + x_m^2)^{\frac{\lambda}{2}}
$$
their action on a test function $\phi \in \mathcal{S}(\mR^m)$ being given by
$$
\la T_\lambda, \phi \ra = \sigma_m \la {\rm Fp} \; r^\mu_+, \Sigma^{(0)}[\phi] \ra 
$$
with $\mu = \lambda +m-1$. In the above expressions ${\rm Fp}\; r^\mu_+$ stands for the classical "finite part" distribution on the real $r$-axis  and $\Sigma^{(0)}$ is the scalar valued generalized spherical mean, defined on scalar valued test functions $\phi(\ux)$ by
$$
\Sigma^{(0)}[\phi] = \frac{1}{\sigma_m} \int_{S^{m-1}} \phi(\ux) \, dS(\uom)
$$
This family $\mathcal{T}$ contains, amongst other ones, the fundamental solutions of the natural powers of the Laplace operator in Euclidean space of odd dimension. As convolution operators they give rise to the traditional Riesz potentials (see e.g. \cite{helgason}). The second family $\mathcal{U} = \{ U_\lambda : \lambda \in \mC\}$ of distributions arises in a natural way by the action of the  Dirac operator $\pux$ on $\mathcal{T}$. The $U_{\lambda}$--distributions thus are typical Clifford analysis constructs: they are Clifford--vector valued, and they also arise as products of $T_{\lambda}$--distributions with the distribution $\uom = \frac{\ux}{|\ux|}$, mentioned above. The action of $U_\lambda$ on a test function $\phi \in \mathcal{S}(\mR^m)$ is given by
$$
\la U_\lambda, \phi \ra  = \sigma_m \la {\rm Fp} \; r^\mu_+, \Sigma^{(1)}[\phi] \ra 
$$
with $\mu = \lambda +m-1$, and where the Clifford--vector valued generalized spherical mean $\Sigma^{(1)}$ is defined on scalar valued test functions $\phi(\ux)$ by
$$
\Sigma^{(1)}[\phi] = \frac{1}{\sigma_m} \int_{S^{m-1}} \uom \ \phi(\ux) \, dS(\uom) 
$$
Typical examples in the $\mathcal{U}$--family are the fundamental solutions of the Dirac operator and of its odd natural powers in Euclidean space of odd dimension.\\[-2mm]

The normalized distributions $T^{*}_\lambda$ and $U^{*}_\lambda$ arise when removing the singularities of $T_\lambda$ and $U_\lambda$  by dividing them by an appropriate Gamma-function showing the same simple poles. These normalized distributions are holomorphic mappings from $\lambda \in \mC$ to the space $\mathcal{S}'(\mR^m)$ of tempered distributions. The scalar $T^{*}_\lambda$ distributions are defined by
\begin{equation}
\left \{
\begin{array}{ll}
\displaystyle{T_\lambda^* = \pi^{\frac{\lambda+m}{2}} \frac{T_\lambda}{\Gamma \left ( \frac{\lambda+m}{2} \right )}}, & \lambda \ne -m-2l\\[5mm]
\displaystyle{T_{-m-2l}^* = \frac{\pi^{\frac{m}{2}-l}}{2^{2l} \Gamma \left ( \frac{m}{2} + l \right )} (-\Delta_m)^l \delta (\ux)}, & l \in \mN_0
\label{defTstar}
\end{array}
\right . 
\end{equation}
while the Clifford--vector valued distributions $U^{*}_\lambda$ are defined by
\begin{equation}
\left \{
\begin{array}{ll}
\displaystyle{U_\lambda^* = \pi^{\frac{\lambda+m+1}{2}} \, \frac{U_\lambda}{\Gamma \left ( \frac{\lambda + m + 1}{2} \right )}}, & \lambda \ne -m-2l-1\\[5mm]
\displaystyle{U_{-m-2l-1}^* = - \frac{\pi^{\frac{m}{2}-l}}{2^{2l+1} \, \Gamma \left ( \frac{m}{2} + l + 1 \right )} \; \pux^{2l+1} \delta(\ux)}, & l \in \mN_0
\label{defUstar}
\end{array}
\right . 
\end{equation}

In this paper we shall also be concerned with the distributions $\pux^k\delta$ and $\pux^kH, k \in \mZ$, where $\delta(\ux)$ stands for the Dirac delta--distribution in $\mR^m$ and $H(\ux)$ for the Hilbert kernel in $\mR^m$ which, through convolution, gives rise to the multidimensional Hilbert transform in the context of Clifford analysis (see e.g. \cite{gilmur}).\\
Let us first introduce the integer powers of the Dirac operator $\pux$.\\
The complex power of the Dirac operator $\pux$ was already introduced in \cite{dss} and further studied in \cite{distrib}. It is a convolution operator defined by
\begin{eqnarray}
	\pux^{\mu}[ \, . \, ]  = \pux^{\mu}\delta \ast [ \, . \, ] & = & \left[\frac{1+e^{i\pi\mu}}{2} \, \frac{2^{\mu}\Gamma\left(\frac{m+\mu}{2}\right)}{\pi^{\frac{m-\mu}{2}}}\; T^{*}_{-m-\mu} - \frac{1-e^{i\pi\mu}}{2} \, \frac{2^{\mu}\Gamma\left(\frac{m+\mu+1}{2}\right)}{\pi^{\frac{m-\mu+1}{2}}}\; U^{*}_{-m-\mu}\right] \ast [ \, . \, ] \nonumber \\[2mm]
	& = & \frac{2^{\mu}}{\pi^{\frac{m}{2}}} \; \mbox{Fp} \, \frac{1}{\vert\ux\vert^{\mu+m}}\left[\frac{1+e^{i\pi\mu}}{2} \, \frac{\Gamma\left(\frac{m+\mu}{2}\right)}{\Gamma\left(-\frac{\mu}{2}\right)}-\frac{1-e^{i\pi\mu}}{2} \, \frac{\Gamma\left(\frac{m+\mu+1}{2}\right)}{\Gamma\left(-\frac{\mu-1}{2}\right)}\uom\right] \ast [ \, . \, ]
\label{defdiracmu}
\end{eqnarray}
	
In particular for natural values of the  parameter $\mu$, the convolution kernel  $\pux^{\mu}\delta$ is given by

\begin{equation}
\left \{ \begin{array}{rcl}
\pux^{2k}\delta & =  & \phantom{-} \displaystyle\frac{2^{2k}\Gamma\left(\frac{m+2k}{2}\right)}{\pi^{\frac{m-2k}{2}}}\; T^{*}_{-m-2k}\\[5mm]
\pux^{2k+1}\delta & = &  - \displaystyle\frac{2^{2k+1}\Gamma\left(\frac{m+2k+2}{2}\right)}{\pi^{\frac{m-2k}{2}}}\; U^{*}_{-m-2k-1}
\end{array} \right . 
\label{natpowdirac}
\end{equation}
which are in accordance with the definitions (\ref{defTstar}) and (\ref{defUstar}). One would be tempted to define for the negative integer powers of the Dirac operator:

\begin{equation}
\left \{ \begin{array}{rcl}
\pux^{-2k}\delta & =  & \phantom{-} \displaystyle\frac{2^{-2k}\Gamma\left(\frac{m-2k}{2}\right)}{\pi^{\frac{m+2k}{2}}}\; T^{*}_{-m+2k}\\[5mm]
\pux^{-2k+1}\delta & = &  - \displaystyle\frac{2^{-2k+1}\Gamma\left(\frac{m-2k+2}{2}\right)}{\pi^{\frac{m+2k}{2}}}\; U^{*}_{-m+2k-1}
\end{array} \right . 
\label{intpowdirac}
\end{equation}
which indeed is a valid definition provided the dimension $m$ is odd. However, if the dimension $m$ is even, the expressions (\ref{intpowdirac}) are not valid for $2k=m, m+2, m+4,\ldots$ in the case of $\pux^{-2k}\delta$ and for $2k=m+2, m+4, \ldots$ in the case of $\pux^{-2k+1}\delta$. For those exceptional parameter values, $\pux^\mu$ is defined as follows:
\begin{equation}
\pux^{-m-n} \delta = E_{m+n}
\label{exepintpowdirac}
\end{equation}
where $E_{m+n}$ is the fundamental solution of the operator $\pux^{m+n}$. For the explicit expression of those fundamental solutions we recall a result of \cite{distrib}.

\begin{proposition}
\label{fusoldirac}
If the dimension $m$ is even, for $n=0,1,2,\ldots$, the fundamental solution $E_{m+n}$ of the operator $\pux^{m+n}$ is given by
$$
\left \{ \begin{array}{rcl}
E_{m+2j} & =  &  (p_{2j} \ln{r} + q_{2j}) \; T^{*}_{2j}\\[1mm]
E_{m+2j+1} & = & (p_{2j+1} \ln{r} + q_{2j+1}) \; U^{*}_{2j+1}
\end{array} \right . \qquad j=0,1,2,\ldots
$$
where the constants $p_n$ and $q_n$ satisfy the recurrence relations
$$
\left \{ \begin{array}{rcl}
p_{2j+2} & =  &  \displaystyle\frac{1}{2j+2} \, p_{2j+1}\\[4mm]
q_{2j+2} & = & \displaystyle\frac{1}{2j+2} \, (q_{2j+1} - \displaystyle\frac{1}{2j+2} \, p_{2j+1})
\end{array} \right . \qquad j=0,1,2,\ldots 
$$
and
$$ 
\left \{ \begin{array}{rcl}
p_{2j+1} & =  & - \displaystyle\frac{1}{2\pi} \, p_{2j}\\[4mm]
q_{2j+1} & = & - \displaystyle\frac{1}{2\pi} \, (q_{2j} - \displaystyle\frac{1}{m+2j} \, p_{2j})
\end{array} \right . \qquad j=0,1,2,\ldots 
$$
with starting values $p_{0} =  - \displaystyle\frac{1}{2^{m-1}\pi^m}$ and $q_{0} = 0$.
\end{proposition}

Next we define the convolution operator $^\mu \mcH$ by
$$
^\mu \mcH [ \, . \, ] = \pux^\mu H \ast [ \, . \, ]
$$
where the convolution kernel $\pux^\mu H$ is given by
$$
\pux^{\mu} H  =  \frac{1-e^{i\pi\mu}}{2} \, \frac{2^{\mu}\Gamma\left(\frac{m+\mu}{2}\right)}{\pi^{\frac{m-\mu}{2}}} \; T^{*}_{-m-\mu} - \frac{1+e^{i\pi\mu}}{2} \, \frac{2^{\mu}\Gamma\left(\frac{m+\mu+1}{2}\right)}{\pi^{\frac{m-\mu+1}{2}}} \; U^{*}_{-m-\mu}
$$
The notation for this new kernel is motivated by the fact that, as shown by a straightforward calculation, it may indeed be obtained as $\pux^{\mu} H  = \pux^{\mu} \delta \ast H$. In particular for natural values of the parameter $\mu$, the convolution kernel $\pux^{\mu} H $ reduces to
\begin{equation}
\left \{ \begin{array}{rcl}
\pux^{2k}H & =  & - 2^{2k}  \displaystyle\frac{\Gamma(\frac{m+2k+1}{2})}{\pi^{\frac{m-2k+1}{2}}} \; U^{*}_{-m-2k} \\[5mm]
\pux^{2k+1}H & = &  2^{2k+1}  \displaystyle\frac{\Gamma(\frac{m+2k+1}{2})}{\pi^{\frac{m-2k-1}{2}}} \; T^{*}_{-m-2k-1} 
\end{array} \right .
\label{natpowhilbertdirac}
\end{equation}
Note that for $\mu=0$ the operator $^0 \mcH$ reduces to the Hilbert transform (see e.g. \cite{gilmur})
$$
^0 \mcH[ \, . \, ] =  \mcH[ \, . \, ] = H \ast [ \, . \, ] = - \frac{2}{\sigma_{m+1}} \, {\rm Fp} \, \frac{\uom}{r^m} \ast  [ \, . \, ]
$$
while for $\mu=1$ the so--called Hilbert--Dirac operator (see e.g.  \cite{fbhds}) is obtained:
$$
^1 \mcH[ \, . \, ] = (-\Delta_m)^\onehalf [ \, . \, ] = \pux H \ast [ \, . \, ] = 2 \frac{\Gamma(\frac{m+1}{2})}{\pi^{\frac{m-1}{2}}} \; T^{*}_{-m-1} \ast [ \, . \, ]
$$

For negative integer parameter values we have
\begin{equation}
\left \{ \begin{array}{rcl}
\pux^{-2k}H & =  & - 2^{-2k}  \displaystyle\frac{\Gamma(\frac{m-2k+1}{2})}{\pi^{\frac{m+2k+1}{2}}} \; U^{*}_{-m+2k} \\[5mm]
\pux^{-2k+1}H & = &  2^{-2k+1}  \displaystyle\frac{\Gamma(\frac{m-2k+1}{2})}{\pi^{\frac{m+2k-1}{2}}} \; T^{*}_{-m+2k-1} 
\end{array} \right .
\label{intpowhilbertdirac}
\end{equation}
which is valid for all natural values of $k$ on condition that the dimension $m$ is even. If $m$ is odd, then the expressions (\ref{intpowhilbertdirac}) fail for $k = \frac{m+1}{2}, \frac{m+1}{2} +1, \ldots$. For these exceptional values we have
$$
\pux^{-m-n} H = F_{m+n}
$$
where $F_{m+n}$ is the fundamental solution of the convolution operator $^{m+n}\mcH = \pux^{m+n}H \ast [\cdot]$, the explicit expression for which is given by the following proposition from \cite{bdbds2}.

\begin{proposition}
\label{fusolhilbert}
If the dimension $m$ is odd, then, for $n=0,1,2,\ldots$, the fundamental solution of $^{m+n}\mcH$ is given by
$$
\left \{ \begin{array}{rcl}
F_{m+2j} & =  &  (p_{2j} \ln{r} + q_{2j})  \; T^{*}_{2j} \\[2mm]
F_{m+2j+1} & = & ( p_{2j+1} \ln{r} + q_{2j+1})    \; U^{*}_{2j+1} 
\end{array} \right . \qquad j=0,1,2,\ldots
$$
with the same constants $(p_n,q_n)$ as in Proposition \ref{fusoldirac}
\end{proposition}



\section{Harmonic and monogenic potentials in $\mR^{m+1}$}
\label{potentials}


In \cite{bdbds1,bdbds2} harmonic and monogenic potentials were studied in upper half--space $\mR_+^{m+1}$ and their distributional boundary values were determined when approaching the boundary $\mR^m$ from that upper half--space. In this section we will consider those potentials also in lower half space $\mR_-^{m+1}$ and investigate the possibility to extend their definition domain across the boundary $\mR^m$.\\[-2mm]

\subsection{The Cauchy kernel}

The Cauchy kernel of Clifford analysis, i.e. the fundamental solution of the generalized Cauchy--Riemann operator $D$:
$$ 
C_{-1}(x_0,\ux) = \onehalf A_{-1}(x_0,\ux) + \onehalf \overline{e_0} B_{-1}(x_0,\ux) = \frac{1}{\sigma_{m+1}} \, \frac{x \overline{e_0}}{|x|^{m+1}} =  \frac{1}{\sigma_{m+1}} \, \frac{x_0 - 
\overline{e_0} \ux}{|x|^{m+1}}
$$
is monogenic in $\mR^{m+1} \setminus \{O\}$, and its two components, which are nothing else but
the traditional Poisson kernels:
\begin{eqnarray*}
A_{-1}(x_0,\ux) & = & P(x_0,\ux) \ = \ \phantom{-} \frac{2}{\sigma_{m+1}} \, \frac{x_0}{|x|^{m+1}} \label{A-1}\\
B_{-1}(x_0,\ux) & = & Q(x_0,\ux) \ = \ - \frac{2}{\sigma_{m+1}} \, \frac{\ux}{|x|^{m+1}} \label{B-1}
\end{eqnarray*}
are conjugate harmonic in the same region $\mR^{m+1} \setminus \{O\}$. For the notion of higher dimensional harmonic conjugate we refer to \cite{fbrdfs} .\\

For $x_0 \neq 0$ we have
$$
\lim_{\ux \rightarrow 0} A_{-1}(x_0,\ux) \ = \  \frac{2}{\sigma_{m+1}} \, \frac{x_0}{|x_0|^{m+1}}
$$
while for $\ux \neq 0$ we have
$$
\lim_{x_0 \rightarrow 0} A_{-1}(x_0,\ux) \ = \ 0
$$
This is in accordance with the distributional boundary values
$$
\left \{ \begin{array}{rcl}
a_{-1}^+(\ux) & = & \lim_{x_0 \rightarrow 0+} A_{-1}(x_0,\ux)  =  \phantom{-} \delta(\ux)\\[2mm]
a_{-1}^-(\ux) & = & \lim_{x_0 \rightarrow 0-} A_{-1}(x_0,\ux)  =  - \delta(\ux) = - a_{-1}^+(\ux)
\end{array} \right .
$$
Note that we will use systematically the superscript notation $\pm$ for denoting the distributional boundary value when approaching $\mR^m$ from the upper, respectively lower, half--space.\\

On the other hand we have, for $x_0 \neq 0$,
$$
\lim_{\ux \rightarrow 0} \, B_{-1}(x_0,\ux) = 0
$$
while for $\ux \neq 0$ there holds
$$
\lim_{x_0 \rightarrow 0} \, B_{-1}(x_0,\ux) = - \frac{2}{\sigma_{m+1}} \, \frac{\uom}{r^{m}}
$$
and, in distributional sense,
$$
\left \{ \begin{array}{rcl}
b_{-1}^+(\ux) & = & \lim_{x_0 \rightarrow 0+} B_{-1}(x_0,\ux)  =  - \frac{2}{\sigma_{m+1}} \, \mbox{Pv} \frac{\uom}{r^{m}} = H(\ux)\\[3mm]
b_{-1}^-(\ux) & = & \lim_{x_0 \rightarrow 0-} B_{-1}(x_0,\ux)  = H(\ux)  = b_{-1}^+(\ux)
\end{array} \right .
$$

\subsection{The downstream potentials}

The downstream potentials are defined recursively by the successive action of the conjugate Cauchy--Riemann operator on the Cauchy kernel $C_{-1}(x_0,\ux)$:

$$
C_{-k-1}(x_0,\ux) =  \onehalf A_{-k-1}(x_0,\ux) + \onehalf \overline{e_0} B_{-k-1}(x_0,\ux) = \overline{D}^k C_{-1}(x_0,\ux), \qquad k=1,2,\ldots
$$

It follows that the downstream potentials $C_{-k-1}$ are monogenic in $\mR^{m+1} \setminus \{ O\}$, and that their components $(A_{-k-1},B_{-k-1})$ are conjugate harmonic in the same region $\mR^{m+1} \setminus \{ O\}$. This is in accordance with the following distributional boundary values ($l=1,2,\ldots$)

$$
\left \{ \begin{array}{rcl}
a_{-2\ell}^+ & = &  (-1)^{\ell-1}  2^{\ell-1} (2\ell-1)!!  \displaystyle\frac{\Gamma \left ( \frac{m+2\ell-1}{2} \right )} {\pi^{\frac{m+1}{2}}} \, {\rm Fp} \displaystyle\frac{1}{r^{m+2\ell-1}}                                \\[4mm]
& = & (-1)^{\ell-1} (2\ell-1)!! (m+1)(m+3)\cdots(m+2\ell-3) \displaystyle\frac{2}{\sigma_{m+1}} \, {\rm Fp} \displaystyle\frac{1}{r^{m+2\ell-1}}   = - \pux^{2\ell-1}H(\ux)  \\[4mm]
a_{-2\ell}^- & = &  a_{-2\ell}^+  \end{array} \right .
$$
$$
\left \{ \begin{array}{rcl}
a_{-2\ell-1}^+ & = & \pux^{2\ell}\delta                            \\[3mm]
a_{-2\ell-1}^- & = &  - \pux^{2\ell}\delta   = -  a_{-2\ell-1}^+  \end{array} \right .
$$
and
$$
\left \{ \begin{array}{rcl}
b_{-2\ell}^+ & = &   - \pux^{2\ell-1} \delta              \\[3mm]
b_{-2\ell}^- & = &  \pux^{2\ell-1} \delta  =  - b_{-2\ell}^+ \end{array} \right .
$$
\vspace*{3mm}
$$
\left \{ \begin{array}{rcl}
b_{-2\ell-1}^+ & =  & (-1)^{\ell-1}  2^{\ell} (2\ell-1)!!  \displaystyle\frac{\Gamma \left ( \frac{m+2\ell+1}{2} \right )}{\pi^\frac{m+1}{2}} \,  {\rm Fp} \displaystyle\frac{\uom}{r^{m+2\ell}}   \\[4mm]
& = & (-1)^{\ell-1} (2\ell-1)!! (m+1)(m+3)\cdots(m+2\ell-1)  \displaystyle\frac{2}{\sigma_{m+1}} \,  {\rm Fp} \displaystyle\frac{\uom}{r^{m+2\ell}}  = \pux^{2\ell} H(\ux) \\[4mm]
b_{-2\ell-1}^-  & = &  b_{-2\ell-1}^+  \end{array} \right .
$$

\subsection{Green's function}

The fundamental solution of the Laplace operator $\Delta_{m+1}$ in $\mR^{m+1}$, sometimes called Green's function, is given by
$$
\frac{1}{2}A_0(x_0,\ux) = - \frac{1}{m-1} \frac{1}{\sigma_{m+1}} \frac{1}{(x_0^2 + r^2)^\frac{m-1}{2}}
$$
Clearly it is a harmonic function in $\mR^{m+1} \setminus \{ O \}$. 
We have for $x_0 \neq 0$
$$
\lim_{\ux \rightarrow 0} \, A_0(x_0,\ux) = - \frac{2}{m-1} \frac{1}{\sigma_{m+1}} \frac{1}{|x_0|^{m-1}}
$$
while for $\ux \neq 0$ we have 
$$
\lim_{x_0 \rightarrow 0+} A_0(x_0,\ux) = \lim_{x_0 \rightarrow 0-} A_0(x_0,\ux) = - \frac{2}{m-1} \frac{1}{\sigma_{m+1}} \frac{1}{r^{m-1}}
$$
This is in accordance with the distributional limits
$$
a_0(\ux)^+ = a_0(\ux)^- = - \frac{2}{m-1} \frac{1}{\sigma_{m+1}} \frac{1}{r^{m-1}}
$$

In the two half--spaces $\mR^{m+1}_+$ and $\mR^{m+1}_-$ separately, a conjugate harmonic to $A_0(x_0,\ux)$ is given by
\begin{equation}
B_0(x_0,\ux) = \frac{2}{\sigma_{m+1}} \, \frac{\ux}{r^m} \, F_m \left ( \frac{r}{x_0} \right )
\label{B0}
\end{equation}
where
$$
F_m(v) = \int_0^v \frac{\eta^{m-1}}{(1+\eta^2)^\frac{m+1}{2}} \, d\eta = \frac{v^m}{m} \,  _2F_1 \left ( \frac{m}{2},\frac{m+1}{2};\frac{m}{2}+1;-v^2 \right )
$$
with $_2F_1$ a standard hypergeometric function (see e.g. \cite{grad}). Note the specific values
$$
F_m(0) = 0
$$
and
$$
F_m(+\infty)= \frac{\sqrt{\pi}}{2} \, \frac{\Gamma(\frac{m}{2})}{\Gamma(\frac{m+1}{2})}
$$
Also note that in low dimensions ($m=2,3$), this function $F_m$ may be expressed in terms of elementary functions (see \cite{bdbds3}).\\

We have for $x_0 \neq 0$
$$
\lim_{\ux \rightarrow 0} \, B_0(x_0,\ux) = 0
$$
while for $\ux \neq 0$ we have to distinguish between
$$
\lim_{x_0 \rightarrow 0+} B_0(x_0,\ux) = \frac{1}{\sigma_m} \frac{\ux}{r^m}
$$
and
$$
\lim_{x_0 \rightarrow 0-} B_0(x_0,\ux) = (-1)^m \, \frac{1}{\sigma_m} \frac{\ux}{r^m}
$$
In distributional sense we also have
$$
\left \{ \begin{array}{rcl}
b_{0}^+ & =  &   \dps\frac{1}{\sigma_m} \dps\frac{\uom}{r^{m-1}} \\[4mm]
b_{0}^-  & = & (-1)^m \,  \dps\frac{1}{\sigma_m} \dps\frac{\uom}{r^{m-1}} = (-1)^m \, b_{0}^+
  \end{array} \right .
$$
This means that the boundary values of $B_0(x_0,\ux)$  from the lower half--space $\mR^{m+1}_-$ depend upon the parity of the dimension considered. It follows that if the dimension $m$ is {\em even}, then $B_0(x_0,\ux)$ can be continued over the boundary $\mR^m$ to a conjugate harmonic function to Green's function $A_0(x_0,\ux)$ in $\mR^{m+1}\setminus \{ O\}$, while if the dimension $m$ is {\em odd} this is not possible and the potential 
$
C_0(x_0,\ux) = \onehalf A_0(x_0,\ux) + \onehalf \overline{e_0} \, B_0(x_0,\ux)
$
remains monogenic in the two half--spaces $\mR^{m+1}_+$ and $\mR^{m+1}_-$ separately.

\begin{remark}
In the upper and lower half of the complex plane the function $\ln(z)$ is a holomorphic potential (or primitive) of the Cauchy kernel $\frac{1}{z}$. By similarity we could say that $C_0(x_0,\ux) = \frac{1}{2} A_0(x_0,\ux) + \frac{1}{2} \overline{e_0} B_0(x_0,\ux)$, being a monogenic potential of the Cauchy kernel $C_{-1}(x_0,\ux)$ and the sum of the fundamental solution $A_0(x_0,\ux)$ of the Laplace operator and its conjugate harmonic $\overline{e_0} B_0(x_0,\ux)$, is a {\em monogenic logarithmic function} in the upper and lower half--spaces $\mR^{m+1}_+$ and $\mR^{m+1}_-$. If $m$ is even then it can even be continued through the boundary $\mR^m$ to a monogenic function in $\mR^{m+1}\setminus \{ O\}$.
\end{remark}

\subsection{The upstream potential $C_1(x_0,\ux)$}

For $m>2$ the upstream potential $A_1(x_0,\ux)$ is given by
$$
A_1(x_0,\ux) = \displaystyle\frac{2}{m-1} \, \displaystyle\frac{1}{\sigma_{m+1}} \, \displaystyle\frac{1}{r^{m-2}} \, F_{m-2} \left ( \displaystyle\frac{r}{x_0} \right )
$$
For $x_0 \neq 0$ we have
$$
\lim_{\ux \rightarrow 0} \, A_1(x_0,\ux) = 0
$$
while, in distributional sense
$$
\left \{ \begin{array}{rcl}
a_{1}^+(\ux) & =  &  \dps\frac{1}{\sigma_m} \, \dps\frac{1}{m-2} \, \dps\frac{1}{r^{m-2}} \\[4mm]
a_{1}^-(\ux)  & = & (-1)^m \, \dps\frac{1}{\sigma_m} \, \dps\frac{1}{m-2} \, \dps\frac{1}{r^{m-2}} = (-1)^m \, a_{1}^+
  \end{array} \right .
$$
So, again, the parity of the dimension $m$ plays a role. If $m$ is {\em even}, then $A_1(x_0,\ux)$ can be continued through the boundary $\mR^m$, except for the origin, to obtain a harmonic function in $\mR^{m+1} \setminus \{ O \}$. On the contrary, when $m$ is {\em odd}, then  $A_1(x_0,\ux)$ is harmonic in both half--spaces $\mR^{m+1}_+$ and $\mR^{m+1}_-$ separately, and there is no way to extend it to a function harmonic in a region crossing the boundary $\mR^m$.\\

For the conjugate harmonic $B_1(0,\ux)$ we have the expression
$$
B_1(x_0,\ux) = \displaystyle\frac{2}{\sigma_{m+1}} \, \displaystyle\frac{x_0 \ux}{r^m} \, F_m \left ( \displaystyle\frac{r}{x_0} \right ) - \displaystyle\frac{2}{\sigma_{m+1}} \, \displaystyle\frac{1}{m-1} \, \displaystyle\frac{\ux}{r^{m-1}}
$$
which clearly is harmonic in $\mR^{m+1} \setminus \{ O \}$.\\
For $x_0 \neq 0$ we have
$$
\lim_{\ux \rightarrow 0} \, B_1(x_0,\ux) = 0
$$
while, in distributional sense
$$
\left \{ \begin{array}{rcl}
b_{1}^+(\ux) & =  &  -  \dps\frac{2}{\sigma_{m+1}} \, \dps\frac{1}{m-1} \, \dps\frac{\uom}{r^{m-2}} \\[4mm]
b_{1}^-(\ux)  & = & b_1^+(\ux)
  \end{array} \right .
$$
In conclusion, we have found that if the dimension $m$ is even ($m>2$), then $A_1(x_0,\ux)$ and $B_1(x_0,\ux)$ are conjugate harmonic in $\mR^{m+1} \setminus \{ O \}$, and 
$
C_1(x_0,\ux) = \onehalf A_1(x_0,\ux) + \onehalf \overline{e_0} \, B_1(x_0,\ux)
$
will be monogenic in the same region $\mR^{m+1} \setminus \{ O \}$. If, on the contrary, the dimension $m$ is odd, then the conjugate harmonicity of $A_1(x_0,\ux)$ and $B_1(x_0,\ux)$ and the monogenicity of $C_1(x_0,\ux)$ only hold in both half--spaces $\mR^{m+1}_+$ and $\mR^{m+1}_-$ separately.

\subsection{The upstream potential $C_2(x_0,\ux)$}

For $m>3$ the upstream potential $A_2(x_0,\ux)$ is given by

\begin{eqnarray*}
A_2(x_0,\ux) & = & \frac{2}{m-1} \frac{1}{\sigma_{m+1}} \frac{x_0}{r^{m-2}} \, F_{m-2} \left ( \frac{r}{x_0} \right ) - \frac{2}{m-1} \frac{1}{m-3} \frac{1}{\sigma_{m+1}} \frac{1}{(x_0^2+r^2)^{\frac{m-3}{2}}}\\
& = & \frac{2}{m-1} \frac{1}{m-2} \frac{1}{\sigma_{m+1}}  \frac{1}{x_0^{m-3}} \, _2F_1\left(\frac{m}{2}-1, \frac{m-1}{2}; \frac{m}{2}; - \frac{r^2}{x_0^2}\right) - \frac{2}{m-1} \frac{1}{m-3} \frac{1}{\sigma_{m+1}} \frac{1}{(x_0^2+r^2)^{\frac{m-3}{2}}}
\end{eqnarray*}

For $x_0 \neq 0$ we have
$$
\lim_{\ux \rightarrow 0} A_2(x_0,\ux) = - \frac{2}{m-1} \frac{1}{m-3} \frac{1}{\sigma_{m+1}} \frac{1}{x_0^{m-3}}
$$
while for $\ux \neq 0$ we have
$$
\lim_{x_0 \rightarrow 0} A_2(x_0,\ux) = - \frac{2}{m-1} \frac{1}{m-3} \frac{1}{\sigma_{m+1}} \frac{1}{r^{m-3}}
$$
This means that $A_2(x_0,\ux)$ is a harmonic function in $\mR^{m+1} \setminus \{ O \}$. Moreover, in distributional sense we have
$$
\left \{ \begin{array}{rcl}
a_{2}^+(\ux) & =  &  -  \dpsf{2}{m-1} \dpsf{1}{m-3} \dps\frac{1}{\sigma_{m+1}} \, {\rm Fp}\dps\frac{1}{r^{m-3}} \\[4mm]
a_{2}^-(\ux)  & = & a_{2}^+(\ux)
  \end{array} \right .
$$
For a conjugate harmonic to $A_2(x_0,\ux)$ we have the following expression:

\begin{eqnarray*}
B_2(x_0,\ux) & = & \frac{1}{\sigma_{m+1}} \frac{\ux (x_0^2+r^2)}{r^{m}} \, F_{m} \left ( \frac{r}{x_0} \right ) - \frac{m-3}{m-1} \frac{1}{\sigma_{m+1}} \frac{\ux}{r^{m-2}} \,   F_{m-2} \left ( \frac{r}{x_0} \right )
\end{eqnarray*}

For $x_0 \neq 0$ we have
$$
\lim_{\ux \rightarrow 0} B_2(x_0,\ux) = 0
$$
while for $\ux \neq 0$ we have to distinguish between
$$
\lim_{x_0 \rightarrow 0+} B_2(x_0,\ux) = \onehalf \frac{1}{\sigma_m} \frac{1}{m-2} \frac{\ux}{r^{m-2}}
$$
and
$$
\lim_{x_0 \rightarrow 0-} B_2(x_0,\ux) = (-1)^m \onehalf \frac{1}{\sigma_m} \frac{1}{m-2} \frac{\ux}{r^{m-2}}
$$

In distributional sense we have
$$
\left \{ \begin{array}{rcl}
b_{2}^+(\ux) & =  &  \dpsf{1}{2}  \dps\frac{1}{m-2}  \dps\frac{1}{\sigma_{m}} \, {\rm Fp}\dps\frac{\uom}{r^{m-3}} \\[4mm]
b_{2}^-(\ux)  & = & (-1)^m \, b_2^+(\ux)
  \end{array} \right .
$$

If the dimension $m$  is even ($m>3$), then $B_2(x_0,\ux)$ becomes harmonic in $\mR^{m+1} \setminus \{ O \}$ and a conjugate harmonic to $A_2(x_0,\ux)$ in the same region $\mR^{m+1} \setminus \{ O \}$, entailing the monogenicity of 
$$
C_2(x_0,\ux) = \frac{1}{2} A_2(x_0,\ux) + \frac{1}{2} \overline{e_0} B_2(x_0,\ux)
$$
in the same region $\mR^{m+1} \setminus \{ O \}$ too.\\
If $m$ is odd ($m>3$), then $B_2(x_0,\ux)$ is a conjugate harmonic to $A_2(x_0,\ux)$ in the half--spaces $\mR^{m+1}_+$ and $\mR^{m+1}_+$ separately, and $C_2(x_0,\ux) $ is monogenic in both half--spaces too.

\subsection{The upstream potentials $C_{k}(x_0,\ux), k=3,\ldots$}

We put, for general $k=1,2,3,\ldots$
$$
A_k(x_0,\ux) = A_k^{\pm}(x_0,\ux) \quad {\rm and} \quad B_k(x_0,\ux) = B_k^{\pm}(x_0,\ux), \quad \ux \in \mR^{m+1}_{\pm}
$$
and define the functions $A_k^{\pm}$ and $ B_k^{\pm}$ in $\mR^{m+1}_{\pm}$ recursively, convolution being taken in the variable $\ux \in \mR^{m}$, by
\begin{eqnarray*}
A_k^+(x_0,\ux)  & = & a_0^+ \ast A_{k-1}^+ \ = \ a_1^+ \ast A_{k-2}^+ \ = \ \ldots \ = \ a_{k-1}^+ \ast A_0^+ \\
& = & b_0^+ \ast B_{k-1}^+ \ = \ b_1^+ \ast B_{k-2}^+ \ = \ \ldots \ = \ b_{k-1}^+ \ast B_0^+ \\
A_k^-(x_0,\ux)  & = & (-1)^m a_0^- \ast A_{k-1}^- \ = \ (-1)^m a_1^- \ast A_{k-2}^- \ = \ \ldots \ = \ (-1)^m a_{k-1}^- \ast A_0^- \\
& = & (-1)^m b_0^- \ast B_{k-1}^- \ = \ (-1)^m b_1^- \ast B_{k-2}^- \ = \ \ldots \ = \ (-1)^m b_{k-1}^- \ast B_0^- \\
B_k^+(x_0,\ux) & = & a_0^+ \ast B_{k-1}^+ \ = \ a_1^+ \ast B_{k-2}^+ \ = \ \ldots \ = \ a_{k-1}^+ \ast B_0^+ \\
& = & b_0^+ \ast A_{k-1}^+ \ = \ b_1^+ \ast A_{k-2}^+ \ = \ \ldots \ = \ b_{k-1}^+ \ast A_0^+ \\
B_k^-(x_0,\ux) & = & (-1)^m a_0^- \ast B_{k-1}^- \ = \ (-1)^m a_1^- \ast B_{k-2}^- \ = \ \ldots \ = \ (-1)^m a_{k-1}^- \ast B_0^- \\
& = & (-1)^m b_0^- \ast A_{k-1}^- \ = \ (-1)^m b_1^- \ast A_{k-2}^- \ = \ \ldots \ = \ (-1)^m b_{k-1}^- \ast A_0^- \\
\end{eqnarray*}
We also put
$$
C_k^{\pm}(x_0,\ux) = \frac{1}{2} A_k^{\pm}(x_0,\ux) + \frac{1}{2} \overline{e_0} B_k^{\pm}(x_0,\ux) \quad {\rm and} \quad C_k(x_0,\ux)= C_k^{\pm}(x_0,\ux), \ \ux \in \mR^{m+1}_{\pm}
$$
It may be verified, to start with, that $A_k^{\pm}(x_0,\ux)$ and $B_k^{\pm}(x_0,\ux)$ are conjugate harmonic potentials, and that $C_k^{\pm}(x_0,\ux)$ is a monogenic potential  of $C_{k-1}^{\pm}(x_0,\ux)$ in the respective half--spaces $\mR^{m+1}_+$ and $\mR^{m+1}_-$ separately.\\

Their distributional boundary values at $\mR^m$ are given by the recurrence relations
\begin{eqnarray*}
a_k^+(\ux) & = & a_0^+ \ast a_{k-1}^+ \ = \ a_1^+ \ast a_{k-2}^+ \ = \ \ldots \ = \  a_{k-1}^+ \ast a_0^+ \\
& = & b_0^+ \ast b_{k-1}^+ \ = \ b_1^+ \ast b_{k-2}^+ \ = \ \ldots \ = \  b_{k-1}^+ \ast b_0^+ \\
a_k^-(\ux) & = & (-1)^m a_0^- \ast a_{k-1}^- \ = \ (-1)^m a_1^- \ast a_{k-2}^- \ = \ \ldots \ = \  (-1)^m a_{k-1}^- \ast a_0^- \\
& = & (-1)^m b_0^- \ast b_{k-1}^- \ = \ (-1)^m b_1^- \ast b_{k-2}^- \ = \ \ldots \ = \  (-1)^m b_{k-1}^- \ast b_0^- \\
b_k^+(\ux) & = & a_0^+ \ast b_{k-1}^+ \ = \ a_1^+ \ast b_{k-2}^+ \ = \ \ldots \ = \  a_{k-1}^+ \ast b_0^+ \\
& = & b_0^+ \ast a_{k-1}^+ \ = \ b_1^+ \ast a_{k-2}^+ \ = \ \ldots \ = \  b_{k-1}^+ \ast a_0^+ \\
b_k^-(\ux) & = & (-1)^m a_0^- \ast b_{k-1}^- \ = \ (-1)^m a_1^- \ast b_{k-2}^- \ = \ \ldots \ = \  (-1)^m a_{k-1}^- \ast b_0^- \\
& = & (-1)^m b_0^- \ast a_{k-1}^- \ = \ (-1)^m b_1^- \ast a_{k-2}^- \ = \ \ldots \ = \  (-1)^m b_{k-1}^- \ast a_0^- 
\end{eqnarray*}
for which the following explicit formulae may be deduced:
$$
\left \{ \begin{array}{lcl}
a_{2k}^+ = a_{2k}^-  & =  & -  \displaystyle\frac{1}{2^{2k+1}} \,  \displaystyle\frac{\Gamma(\frac{m-2k-1}{2})}{\pi^{\frac{m+2k+1}{2}}} \; T^{*}_{-m+2k+1} \\[5mm]
& = & - \displaystyle\frac{2}{\sigma_{m+1}} \, \displaystyle\frac{1}{(2k-1)!!} \, \displaystyle\frac{1}{(m-1)(m-3)\cdots(m-2k-1)} \, {\rm Fp} \, r^{-m+2k+1}\\[5mm]
& \phantom{=} & (2k \neq m-1,m+1,m+3,\ldots)\\[5mm]
a_{2k-1}^+ = (-1)^m a_{2k-1}^- & =  & \phantom{-} \displaystyle\frac{1}{2^{2k}} \, \displaystyle\frac{\Gamma(\frac{m-2k}{2})}{\pi^{\frac{m+2k}{2}}} \; T^{*}_{-m+2k}\\[5mm]
& = & \displaystyle\frac{1}{2^{k-1}} \, \displaystyle\frac{1}{\sigma_m} \, \displaystyle\frac{1}{(k-1)!} \, \displaystyle\frac{1}{(m-2)(m-4)\cdots(m-2k)} \, {\rm Fp} \, r^{-m+2k}\\[5mm]
& \phantom{=} & (2k \neq m,m+2,m+4,\ldots)
\end{array} \right .
$$
$$
\left \{ \begin{array}{lcl}
b_{2k}^+ = (-1)^m b_{2k}^-& = & \phantom{-}  \displaystyle\frac{1}{2^{2k+1}}  \, \displaystyle\frac{\Gamma(\frac{m-2k}{2})}{\pi^{\frac{m+2k+2}{2}}} \; U^{*}_{-m+2k+1} \\[5mm]
& \phantom{=} & (2k \neq m,m+2,m+4,\ldots)\\[5mm]
& = & \displaystyle\frac{1}{2^{k}}  \, \displaystyle\frac{1}{\sigma_m}  \, \displaystyle\frac{1}{k!}  \, \displaystyle\frac{1}{(m-2)(m-4)\cdots(m-2k)}  \, \uom \, {\rm Fp} \, r^{-m+2k+1}\\[5mm]
b_{2k-1}^+ = b_{2k-1}^- & = &  -  \displaystyle\frac{1}{2^{2k}} \,  \displaystyle\frac{\Gamma(\frac{m-2k+1}{2})}{\pi^{\frac{m+2k+1}{2}}} \; U^{*}_{-m+2k} \\[5mm]
& = & - \displaystyle\frac{2}{\sigma_{m+1}} \, \displaystyle\frac{1}{(2k-1)!!} \, \displaystyle\frac{1}{(m-1)(m-3)\cdots(m-2k+1)} \, \uom \, {\rm Fp} \, r^{-m+2k}\\[5mm]
& \phantom{=} & (2k \neq m+1,m+3,\ldots)
\end{array} \right .
$$
For the above mentioned exceptional values, which occur the sooner the dimension is lower (see \cite{bdbds3}), the distributional boundary values are given by
$$
\left \{ \begin{array}{lcl}
a_{m+2j-1}^+ = a_{m+2j-1}^-  & =  & -  F_{m+2j} \quad j=0,1,2,\ldots (m \ {\rm odd})\\[3mm]
a_{m+2j-1}^+ = a_{m+2j-1}^- & =  & \phantom{-} E_{m+2j} \quad j=0,1,2,\ldots (m \ {\rm even})
\end{array} \right .
$$
$$
\left \{ \begin{array}{lcl}
b_{m+2j}^+ = b_{m+2j}^- & = & - E_{m+2j+1} \quad j=0,1,2,\ldots (m \ {\rm even})\\[3mm]
b_{m+2j}^+ = b_{m+2j}^- & = & \phantom{-}  F_{m+2j+1} \quad j=0,1,2,\ldots  (m \ {\rm odd})\\
\end{array} \right .
$$
with, see \cite{bdbds2} and Section 2,
$$
E_{m+2j} = F_{m+2j} = (p_{2j} \ln{r} + q_{2j}) \, \displaystyle\frac{\pi^{\frac{m+2j}{2}}}{\Gamma(\frac{m+2j}{2})} \, {\rm Fp} \, r^{2j}
$$
$$
E_{m+2j+1} = F_{m+2j+1} = (p_{2j+1} \ln{r} + q_{2j+1}) \, \displaystyle\frac{\pi^{\frac{m+2j+2}{2}}}{\Gamma(\frac{m+2j+2}{2})} \, \uom \, {\rm Fp} \, r^{2j+1}
$$
where the constants $p_n$ and $q_n$ satisfy the recurrence relations of Propositions \ref{fusoldirac} and \ref{fusolhilbert}.\\

These distributional limits show the following properties.
\begin{lemma}
\label{lem58}
One has for $k=1,2,\ldots$
\begin{itemize}
\item[(i)] $- \pux a_k^+ = b_{k-1}^+$; $- \pux b_k^+ = a_{k-1}^+$
\item[(ii)]$- \pux a_k^- = b_{k-1}^-$; $- \pux b_k^- = a_{k-1}^-$
\item[(iii)] $\mathcal{H} \left [ a_k^+ \right ] = b_{-1}^+ \ast a_k^+ = b_k^+$; $\mathcal{H} \left [ b_k^+ \right ] = b_{-1}^+ \ast b_k^+ = a_k^+$
\item[(iv)] $\mathcal{H} \left [ a_k^- \right ] = b_{-1}^- \ast a_k^- = (-1)^m b_k^-$; $\mathcal{H} \left [ b_k^- \right ] = b_{-1}^- \ast b_k^- = (-1)^m a_k^-$
\end{itemize}
\end{lemma}

In conclusion we can state that if $m$ is even, then for all $k=1,2,\ldots$ the potentials $A_k(x_0,\ux)$ and  $B_k(x_0,\ux)$ are conjugate harmonic, and   $C_k(x_0,\ux)$ is monogenic, in $\mR^{m+1} \setminus \{ O \}$. If $m$ is odd, then  for all $k=1,2,\ldots$ the potentials $C_k(x_0,\ux)$ are monogenic in the half--spaces $\mR^{m+1}_+$ and $\mR^{m+1}_-$ separately, the potentials $A_{2k}(x_0,\ux)$ and  $B_{2k-1}(x_0,\ux)$ being harmonic in $\mR^{m+1} \setminus \{ O \}$, while $A_{2k-1}(x_0,\ux)$ and  $B_{2k}(x_0,\ux)$ are harmonic in both half--spaces separately.



\section{Representation of $\pux^n \delta(\ux)$ and $\pux^n H(\ux), n \in \mN$}


In the previous section we have listed the distributional boundary values in $\mR^m$, both from the upper and from the lower half--space, of the harmonic and monogenic potentials considered. In this section we change the viewpoint and aim at representing those distributions in $\mR^m$ as the difference, sometimes called the {\em jump} over $\mR^m$, of monogenic functions in both half--spaces. 


\subsection{Representation of $\delta(\ux)$ and $H(\ux)$}

From subsection 3.1 we know that, putting
$$
\lim_{x_0 \rightarrow 0+} \, C_{-1}(x_0,\ux) = c_{-1}^+(\ux) = \onehalf a_{-1}^+(\ux) + \onehalf \enb b_{-1}^+(\ux)
$$
and similarly
$$
\lim_{x_0 \rightarrow 0-} \, C_{-1}(x_0,\ux) = c_{-1}^-(\ux) = \onehalf a_{-1}^-(\ux) + \onehalf \enb b_{-1}^-(\ux)
$$
there holds
$$
c_{-1}^+(\ux) - c_{-1}^-(\ux) = a_{-1}^+(\ux) = \delta(\ux)
$$
and
$$
c_{-1}^+(\ux) + c_{-1}^-(\ux) = \enb \, b_{-1}^+(\ux) = \enb \, H(\ux)
$$
In this way the delta distribution $\delta(\ux)$ and the Hilbert kernel $H(\ux)$ in $\mR^m$ are represented as the couple of monogenic functions:
\begin{equation}
\label{repdelta}
\delta(\ux) \longleftrightarrow (C_{-1}^+(x_0,\ux),C_{-1}^-(x_0,\ux) ) = \left( \frac{1}{\sigma_{m+1}} \frac{x_0 - \enb \ux}{|x|^{m+1}}       ,  \frac{1}{\sigma_{m+1}} \frac{x_0 - \enb \ux}{|x|^{m+1}}       \right)
\end{equation}
and 
\begin{equation}
\label{rephilbert}
H(\ux) \longleftrightarrow (e_0 \, C_{-1}^+(x_0,\ux), \enb \, C_{-1}^-(x_0,\ux) ) = \left( \frac{1}{\sigma_{m+1}} \frac{x_0 e_0 -  \ux}{|x|^{m+1}}       ,  \frac{1}{\sigma_{m+1}} \frac{x_0 \enb +  \ux}{|x|^{m+1}}       \right)
\end{equation}
where we recall that 
$$
T(\ux)  \longleftrightarrow (F(x_0,\ux),G(x_0,\ux))
$$
stands for
$$
T(\ux) = \lim_{x_0 \rightarrow 0+} \, F(x_0,\ux) - \lim_{x_0 \rightarrow 0-} \, G(x_0,\ux)
$$
the limits being taken in distributional sense and
the functions $F(x_0,\ux)$ and $G(x_0,\ux)$ being monogenic in the respective half--spaces $\mR^{m+1}_+$ and $\mR^{m+1}_-$.\\

The  above representations \eqref{repdelta}, \eqref{rephilbert} of $\delta(\ux)$ and $H(\ux)$ respectively, are in fact nothing else but a reformulation of the well--known Plemejl--Sokhotski formulae in Clifford analysis; they are the multidimensional counterparts to the hyperfunctions \eqref{hyperdelta} and \eqref{hyperhilbert} on the real line.

\subsection{Representation of $\pux \delta(\ux)$ and $\pux H(\ux)$}

From subsection 3.2 we can directly deduce that, with similar definitions as above for $c_{-2}^+(\ux)$ and $c_{-2}^-(\ux)$,
$$
c_{-2}^+(\ux) - c_{-2}^-(\ux) = \enb \, b_{-2}^+(\ux) = - \enb \, \pux \delta(\ux)
$$
and
$$
c_{-2}^+(\ux) + c_{-2}^-(\ux) = a_{-2}^+(\ux) = - \pux H(\ux)
$$
yielding the representations
$$
\pux \delta(\ux) \longleftrightarrow (\enb \, C_{-2}^+(x_0,\ux), \enb \,C_{-2}^-(x_0,\ux))
$$
and   
$$
\pux H(\ux)  \longleftrightarrow (- C_{-2}^+(x_0,\ux), C_{-2}^-(x_0,\ux))
$$  

There is also an indirect way, using the Plemelj--Sokhotski formulae, to obtain the same representation. We indeed have
$$
\pux \delta(\ux) = (c_{-1}^+(\ux) - c_{-1}^-(\ux)) \ast \pux \delta(\ux) = a_{-1}^+(\ux) \ast \pux \delta(\ux) = \pux  a_{-1}^+(\ux) \ast \delta(\ux) 
$$
$$
= - b_{-2}^+(\ux) \ast \delta(\ux) = \enb \, (c_{-2}^+(\ux) - c_{-2}^-(\ux) ) \ast  \delta(\ux) = \enb \, (c_{-2}^+(\ux) - c_{-2}^-(\ux) )
$$
and similarly
$$
\pux H(\ux) = (c_{-1}^+(\ux) - c_{-1}^-(\ux)) \ast \pux H(\ux) = a_{-1}^+(\ux) \ast \pux H(\ux) = \pux  a_{-1}^+(\ux) \ast H(\ux) 
$$
$$
= - b_{-2}^+(\ux) \ast H(\ux) = - \mcH[b_{-2}^+(\ux)] = - a_{-2}^+(\ux) = -c_{-2}^+(\ux) - c_{-2}^-(\ux)
$$

\begin{remark}
The distribution $\pux H(\ux)$ is special. The Dirac operator $\pux$ and the Hilbert kernel $H$ both being vector--valued, the distribution $\pux H(\ux)$ is, surprisingly, scalar--valued, and,  as already mentioned in Section 2, it is the so--called Hilbert--Dirac kernel which, through convolution, gives rise to the well--known scalar pseudodifferential operator ''square root of the Laplacian'' (see e.g. \cite{fbhds}):
$$
\pux H(\ux) =  H(\ux) \pux = (-\Delta_m)^{\onehalf} \delta = - \frac{2}{\sigma_{m+1}} \, {\rm Fp} \, \frac{1}{r^{m+1}}
$$
for which it indeed holds that
$$
- \Delta_m = (-\Delta_m)^{\onehalf}(-\Delta_m)^{\onehalf}
$$
Also the distribution $\pux\delta(\ux)$ is special since it can be expressed as, see \cite{bdbds2},
$$
\pux\delta(\ux) = (-\Delta_m)^{\onehalf} H
$$
These formulae nicely illustrate the symmetric role played by the $\delta$ and $H$ distributions in $\mR^m$.
\end{remark}

\subsection{Representation of $\pux^n \delta(\ux)$ and $\pux^n H(\ux), n \in \mN$}

From subsection 3.2 it follows that
$$
c_{-2\ell}^+(\ux) - c_{-2\ell}^-(\ux) = \enb \, b_{-2\ell}^+(\ux) = - \enb \, \pux^{2\ell-1} \delta(\ux)
$$
and
$$
c_{-2\ell}^+(\ux) + c_{-2\ell}^-(\ux) = a_{-2\ell}^+(\ux) = - \pux^{2\ell-1} H(\ux)
$$
leading to the representations
$$
\pux^{2\ell-1} \delta(\ux) \longleftrightarrow (\enb \, C_{-2\ell}^+(x_0,\ux), \enb \,C_{-2\ell}^-(x_0,\ux))
$$
and   
$$
\pux^{2\ell-1} H(\ux)  \longleftrightarrow (- C_{-2\ell}^+(x_0,\ux), C_{-2\ell}^-(x_0,\ux))
$$  

It also follows that
$$
c_{-2\ell-1}^+(\ux) - c_{-2\ell-1}^-(\ux) =  a_{-2\ell-1}^+(\ux) =  \pux^{2\ell} \delta(\ux)
$$
and
$$
c_{-2\ell-1}^+(\ux) + c_{-2\ell-1}^-(\ux) = \enb \, b_{-2\ell-1}^+(\ux) =  \enb \, \pux^{2\ell} H(\ux)
$$
leading to the representations
$$
\pux^{2\ell} \delta(\ux) \longleftrightarrow ( C_{-2\ell-1}^+(x_0,\ux), C_{-2\ell-1}^-(x_0,\ux))
$$
and   
$$
\pux^{2\ell} H(\ux)  \longleftrightarrow ( e_0 \, C_{-2\ell-1}^+(x_0,\ux), \enb \, C_{-2\ell-1}^-(x_0,\ux))
$$  

\begin{remark}
Here we have obtained the representation of the scalar distributions $\pux^{2\ell-1} H(\ux)$, which, by convolution, yield the half--integer powers of the Laplace operator:
$$
\pux^{2\ell-1} H(\ux) = (-\Delta_m)^{\ell - \onehalf} \delta(\ux)  \quad , \quad \ell=1,2,\ldots
$$
and of the vector distributions $\pux^{2\ell-1} \delta(\ux)$, which may be expressed in a similar way:
$$
\pux^{2\ell-1} \delta(\ux) = (-\Delta_m)^{\ell - \onehalf} H(\ux)
$$
\end{remark}



\section{Representation of $\pux^{-n} \delta(\ux)$ and $\pux^{-n} H(\ux), n \in \mN$}


Recalling the definitions (\ref{intpowdirac}) and (\ref{exepintpowdirac}) of the negative integer powers of the Dirac operator and comparing them with the distributional boundary values of the upstream potentials obtained in subsection 3.6, it is clear that for $m$ odd and for $m$ even with $2k \neq m,m+2, \ldots$
$$
\pux^{-2k}\delta = a_{2k-1}^+(\ux)  \quad  {\rm and} \quad  \pux^{-2k-1}\delta = - b_{2k}^+(\ux) 
$$
and also for $m$ even and $j=0,1,2,\ldots$
$$
\pux^{-m-2j}\delta = a_{m+2j-1}^+(\ux)  \quad  {\rm and} \quad \pux^{-m-2j-1}\delta = - b_{m+2j}^+(\ux)
$$
In a similar way we find that for $m$ even and for $m$ odd with $2k \neq m+1, m+3, \ldots$
$$
\pux^{-2k}H = b_{2k-1}^+(\ux)  \quad  {\rm and} \quad  \pux^{-2k-1}H = - a_{2k}^+(\ux) 
$$
and also for $m$ odd and $j=0,1,2,\ldots$
$$
\pux^{-m-2j}H = - a_{m+2j-1}^+(\ux)  \quad  {\rm and} \quad \pux^{-m-2j-1}H =  b_{m+2j}^+(\ux)
$$

\subsection{Representation of $\pux^{-1} \delta(\ux)$ and $\pux^{-1} H(\ux)$}

In particular we have
$$
c_0^+(\ux) - c_0^-(\ux) = \onehalf \enb \, (1-(-1)^m) \, b_0^+(\ux)
$$
and
$$
c_0^+(\ux) + c_0^-(\ux) =  a_0^+(\ux) + \onehalf \enb \, (1+(-1)^m) \, b_0^+(\ux)
$$
So, if $m$ is odd, we find
$$
c_0^+(\ux) - c_0^-(\ux) =  \enb \, b_0^+(\ux) = e_0 \, \pux^{-1} \delta(\ux) 
$$
and
$$
c_0^+(\ux) + c_0^-(\ux) =  a_0^+(\ux) = - \pux^{-1}H
$$
and the corresponding representations
\begin{equation}
\label{repminusonedirac}
\pux^{-1} \delta(\ux)  \longleftrightarrow ( \enb \, C_0^+(x_0,\ux), \enb \, C_0^-(x_0,\ux) )
\end{equation}
and
\begin{equation}
\label{repminusonehilbert}
\pux^{-1} H(\ux) \longleftrightarrow (- C_0^+(x_0,\ux),  C_0^-(x_0,\ux) )
\end{equation}
Still under the assumption that $m$ is odd, these representations may also be obtained indirectly by means of the Plemelj--Sokhotski formulae:
$$
\pux^{-1} \delta(\ux) = (c_{-1}^+ - c_{-1}^-) \ast \pux^{-1} \delta = a_{-1}^+ \ast \pux^{-1} \delta = \pux^{-1} a_{-1}^+ \ast \delta = - b_0^+ \ast \delta = - e_0 \, (c_0^+ - c_0^-) \ast \delta = \enb \, (c_0^+ - c_0^-)
$$
and
$$
\pux^{-1} H(\ux) = (c_{-1}^+ - c_{-1}^-) \ast \pux^{-1} H = a_{-1}^+ \ast \pux^{-1} H = \pux^{-1} a_{-1}^+ \ast H = - b_0^+ \ast H = - a_0^+ = - (c_0^+ + c_0^-)
$$

However, if $m$ is even, a representation of this kind, involving the monogenic potential $C_0(x_0,\ux)$, for the distributions $\pux^{-1} \delta(\ux)$ and $\pux^{-1} H(\ux)$ is not possible. In this case we have to restrict ourselves to the mere Plemelj--Sokhotsky representations
$$
\pux^{-1} \delta(\ux)  \longleftrightarrow \left(    (C_{-1}(x_0, \cdot) \ast   \pux^{-1} \delta(\cdot))^+(\ux)  \  ,   \     (C_{-1}(x_0, \cdot) \ast   \pux^{-1} \delta(\cdot))^-(\ux)   \right)
$$
and
$$
\pux^{-1} H(\ux) \longleftrightarrow \left(      (C_{-1}(x_0, \cdot) \ast   \pux^{-1} H(\cdot))^+(\ux)  \  ,   \     (C_{-1}(x_0, \cdot) \ast   \pux^{-1} H(\cdot))^-(\ux)             \right)
$$

\begin{remark}
Taking into account that (i) $\pux^{-1}\delta = - E_{1}$, with $E_1(\ux) = - \frac{1}{\sigma_m} \frac{\ux}{r^m}$ the fundamental solution of the Dirac operator $\pux$ in $\mR^m$, (ii) $C_0(x_0,\ux)$ is the monogenic logarithmic function in the half--spaces $\mR^{m+1}_+$ and $\mR^{m+1}_-$, (iii) the Heaviside function $Y(x)$ can be seen as the fundamental solution of the differential operator $\frac{d}{dx}$ on the real line, it is justified to see the hyperfunction representation \eqref{repminusonedirac} in the case where $m$ is odd, as the multidimensional counterpart to the Heaviside hyperfunction \eqref{hyperheaviside} in the complex plane where $m=1$. 
\end{remark}

\subsection{Representation of $\pux^{-2} \delta(\ux)$ and $\pux^{-2} H(\ux)$}

We have
$$
c_1^+(\ux) - c_1^-(\ux) = \onehalf (1-(-1)^m) \, a_1^+(\ux)
$$
and
$$
c_1^+(\ux) + c_1^-(\ux) =  \onehalf (1+(-1)^m) \, a_1^+(\ux) + \enb \, b_1^+(\ux)
$$
So, if $m$ is odd, we find
$$
c_1^+(\ux) - c_1^-(\ux) =  a_1^+(\ux) = \pux^{-2} \delta(\ux) 
$$
and
$$
c_1^+(\ux) + c_1^-(\ux) =  \enb \, b_1^+(\ux) = \enb \, \pux^{-2}H(\ux)
$$
and the corresponding representations
\begin{equation}
\label{repminustwodirac}
\pux^{-2} \delta(\ux)  \longleftrightarrow ( C_1^+(x_0,\ux), C_1^-(x_0,\ux) )
\end{equation}
and
\begin{equation}
\label{repminustwohilbert}
\pux^{-2}H(\ux) \longleftrightarrow (e_0 \, C_1^+(x_0,\ux), \enb \, C_1^-(x_0,\ux) )
\end{equation}
Similar remarks concerning the case where $m$ is even and the Plemelj--Sokhotsky approach can be made as in the preceding subsection.

\begin{remark}
The Clifford hyperfunctions \eqref{repminustwodirac} and \eqref{repminustwohilbert} are the multidimensional counterparts to the complex hyperfunctions
$$
x Y(-x) \longleftrightarrow (\frac{1}{2\pi i} z (\ln{z} - 1) , \frac{1}{2\pi i} z (\ln{z} - 1) )
$$
and
$$
x (\ln{|x|} -1) \longleftrightarrow (\onehalf z (\ln{z} - 1) , - \onehalf z (\ln{z} - 1) )
$$
\end{remark}

\subsection{Representation of $\pux^{-n} \delta(\ux)$ and $\pux^{-n} H(\ux), n \in \mN$}

In view of the results obtained in the preceding subsections, we assume from the start that $m$ is odd. Then we have
$$
c_{2k-1}^+(\ux) - c_{2k-1}^-(\ux) =  a_{2k-1}^+(\ux) = \pux^{-2k} \delta(\ux)
$$
and
$$
c_{2k-1}^+(\ux) + c_{2k-1}^-(\ux) =  \enb \, b_{2k-1}^+(\ux) = \enb \, \pux^{-2k} H(\ux)
$$
leading to the representations
$$
\pux^{-2k} \delta(\ux)  \longleftrightarrow ( C_{2k-1}^+(x_0,\ux), C_{2k-1}^-(x_0,\ux) )
$$
and
$$
\pux^{-2k}H(\ux) \longleftrightarrow (e_0 \, C_{2k-1}^+(x_0,\ux), \enb \, C_{2k-1}^-(x_0,\ux) )
$$
We also have
$$
c_{2k}^+(\ux) - c_{2k}^-(\ux) =  \enb \, b_{2k}^+(\ux) = \enb \, (-\pux^{-2k-1} \delta(\ux))
$$
and
$$
c_{2k}^+(\ux) + c_{2k}^-(\ux) =  a_{2k}^+(\ux) = - \pux^{-2k-1} H(\ux)
$$
leading to the representations
$$
\pux^{-2k-1} \delta(\ux)  \longleftrightarrow ( \enb \, C_{2k}^+(x_0,\ux), \enb \, C_{2k}^-(x_0,\ux) )
$$
and
$$
\pux^{-2k-1}H(\ux) \longleftrightarrow (- C_{2k}^+(x_0,\ux),  C_{2k}^-(x_0,\ux) )
$$
A similar remark as in the preceding subsections concerning the case where $m$ is odd applies also here.



\end{document}